\documentclass[12pt,twoside]{amsart} 

\usepackage{amssymb,latexsym,bbm,graphicx,epsfig,epic,eepic,oldgerm,psfrag}
\usepackage{amsmath,amsfonts,amscd,mathrsfs,amsthm}
\usepackage{tikz,wasysym} 
\usepackage{marvosym}          
\usepackage{a4wide}

\theoremstyle{plain}

\newtheorem*{remark*}{Remark}
\newtheorem*{remarks*}{Remarks}

\newtheorem*{example*}{Example}
\newtheorem*{examples*}{Examples}

\newtheorem*{definition*}{Definition}

\newtheorem*{EPrinciple*}{Eliashberg's Principle}

\newtheorem{inflation.lemma}[theorem]{Inflation Lemma}

\numberwithin{figure}{section}
\numberwithin{equation}{section}

\def\B{\operatorname{B}}

\newcommand{\ballvier}{\overset%
{\raisebox{-0.7ex}[0.0ex][-0.3ex]{\mbox{$\!\! \scriptscriptstyle\circ$}} \mskip2mu} {\operatorname{B}^{ \raisebox{.0ex}{$\scriptstyle 4$}}}}

\newcommand{\ballzwein}{\overset%
{\raisebox{-0.7ex}[0ex][-0.3ex]{\mbox{$\!\!\!\!\!\!\: \scriptscriptstyle\circ$}} \mskip2mu} {\operatorname{B}^{ \raisebox{.0ex}{$\scriptstyle 2n$}}}}

\newcommand{\Kcirc}{\overset%
{\raisebox{-.3ex}[0ex][-.3ex]{\mbox{$\scriptscriptstyle\circ$}} \mskip-5mu} K}

\def\1{\:\!}
\def\2{\;\!}
\def\s{\smallskip}
\def\m{\medskip}
\def\eps{\varepsilon}

\def\Vol{\operatorname {Vol}\2}

\def\Diffc0{\operatorname{Diff^c_0}}

\def\Hamc{\operatorname{Ham^c}}

\def\area{\operatorname{area}}

\def\PD{\operatorname{PD}}
\def\C{\operatorname{C}}

\def\mm{\mathbf{m}}

\def\ww{\boldsymbol{w}}

\def\Bl{\operatorname{Bl}}

\def\Emb{\operatorname{Emb}}

\def\ga{\alpha}

\def\gg{\gamma}

\def\gve{\varepsilon}
\def\gf{\varphi}

\def\go{\omega}

\def\B{\operatorname{B}}

\def\C{\operatorname{C}}
\def\E{\operatorname{E}}

\def\cc{{\mathcal C}}

\def\cl{{\mathcal L}}

\def\CC{\mathbbm{C}}

\def\NN{\mathbbm{N}}

\def\RR{\mathbbm{R}}

\def\ZZ{\mathbbm{Z}}

\def\CP{\operatorname{\mathbbm{C}P}}

\def\pp{\partial}

\newcommand{\se}{\overset%
{\raisebox{-.2ex}[0ex][-.2ex]{\mbox{$\scriptstyle s$}} \mskip 2.5mu}\hookrightarrow}

\def\ni{\noindent}
\def\b{\bigskip}
\def\m{\medskip}

\def\id{\mbox{id}}

\definecolor{amber}{rgb}{1.0, 0.75, 0.0}\definecolor{amber}{rgb}{1.0, 0.75, 0.0}

\definecolor{Green1}{rgb}{0.0, 0.5, 0.0}
\definecolor{Green2}{rgb}{0.07, 0.53, 0.03}

\begin{document}

\title[]{Dusa Mc\,Duff and symplectic geometry}

\author{Felix Schlenk} \thanks{FS partially supported by SNF grant 200020-144432/1.}
\address{(F.~Schlenk)
Institut de Math\'ematiques,
Universit\'e de Neuch\^atel,
Rue \'Emile Argand~11,
2000 Neuch\^atel,
Switzerland}
\email{schlenk@unine.ch}                                              

\date{\today}


\maketitle

Dusa McDuff has led three mathematical lives.
In her early twenties she worked on von Neumann algebras, and in~\cite{69a, 69b} 
established the existence of uncountably many different algebraic types of ${\rm II}_1$-factors.
After an inspiring six months studying with Gel'fand in Moscow and a 
two-year post-doc in Cambridge studying topology,  
she wrote a ``second thesis" while working with Graeme Segal     
on configuration spaces and the group-completion theorem,  
and investigated the topology of various diffeomorphism groups, 
in particular symplectomorphism groups.  
This brought her to symplectic geometry, which was revolutionized around~1985  
by Gromov's introduction of $J$-holomorphic curves~\cite{Gr85}. 

In this short text I will describe some of McDuff's wonderful contributions to 
symplectic geometry.
After reviewing what is meant by `symplectic'
I will mostly focus on her work on symplectic embedding problems.
Some of her other results in symplectic geometry are discussed at the end.
More personal texts about Dusa can be found in~\cite{11}.
Parts of this text overlap with the `Perspective'
in~\cite{11} written jointly with Leonid Polterovich.

\section{Symplectic}

There are many strands to and from symplectic geometry. 
The most important ones are classical mechanics and algebraic geometry.
I do not list these strands 
but refer you to~\cite{McSal.book1, Sch18}.
Here, I simply give the 

\m \ni
{\bf Definition 1.}
Let $M$ be a smooth manifold. 
A {\it symplectic form} on $M$ is a non-degenerate closed 2-form $\omega$.
A diffeomorphism $\gf$ of $M$ is {\it symplectic} (or a {\it symplectomorphism})
if $\gf^* \omega = \omega$.

\s
The non-degeneracy condition implies that symplectic manifolds are even-dimensional.
An example is $\RR^{2n}$ with the constant differential $2$-form 
$$
\go_0 \,=\, \sum_{i=1}^n dx_i \wedge dy_i .
$$
Other examples are surfaces endowed with an area form, their products, 
and K\"ahler manifolds.

\begin{figure}[h] 
 \begin{center}
  \psfrag{+}{$+$}
  \psfrag{-}{$-$}
  \leavevmode\epsfbox{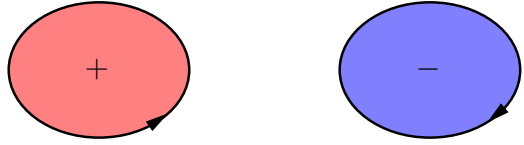} 
 \end{center}
 \caption{The sign of the signed area of an embedded closed curve in $\RR^2$}
 \label{fig.orient}
\end{figure}

\m
If you begin your first lecture on symplectic geometry like this, 
you may very well find yourself alone the following week.
You may thus prefer to start in a more elementary way:
Let $\gg$ be a closed oriented piecewise smooth curve in~$\RR^2$.
If $\gg$ is embedded, assign to~$\gg$ the signed area of the disc $D$ bounded by~$\gg$,
namely \textcolor{red}{$\area (D)$} or \textcolor{blue}{$-\area (D)$},
as in Figure~\ref{fig.orient}. 

%

\begin{figure}[h]
 \begin{center}
  \leavevmode\epsfbox{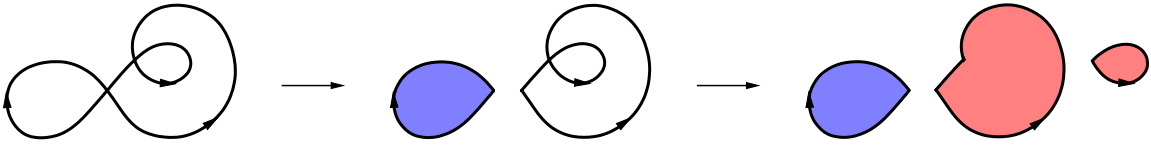}
 \end{center}
 \caption{Splitting a closed curve into embedded pieces}
 \label{fig.split}
\end{figure}

\ni
If $\gg$ is not embedded, successively decompose $\gg$ into closed embedded pieces as illustrated in 
Figure~\ref{fig.split}, 
and define $A (\gg)$ as the sum of the signed areas of these pieces.

%

\m \ni
{\bf Definition 2.}
The standard symplectic structure of~$\RR^{2n}$ is the map  
$$
A (\gg) \,=\, \sum_{i=1}^n A(\gg_i), \quad\; \gg = \left( \gg_1, \dots, \gg_n \right) \subset \RR^{2n} .
$$
A symplectomorphism $\gf$ of $\RR^{2n}$ is a diffeomorphism that preserves the signed area of 
closed curves:
$$
A (\gf (\gg)) \,=\, A(\gg) \quad \mbox{ for all closed curves } \gg \subset \RR^{2n} .
$$
A symplectic structure on a manifold $M$ is an atlas whose transition functions are (local) symplectomorphisms,
and a symplectomorphism of~$M$ is then a diffeomorphism that preserves this local structure.

\m
The standard symplectic structure of~$\RR^{2n}$ is thus given by assigning to 
a closed curve
the sum of the signed areas of 
the $n$~curves obtained by projecting to the coordinate planes~$\RR^2(x_i,y_i)$.
And a symplectic structure on a manifold is a coherent way of assigning a 
signed area to sufficiently local closed curves.
The equivalence of the two definitions follows from 

\m \ni
{\bf Darboux's Theorem.}
{\it Around every point of a symplectic manifold $(M,\go)$ there exists 
a coordinate chart~$\gf$ such that $\gf^* \go_0 = \go$.
}

\m
The group of symplectomorphisms of a symplectic manifold is very large. 
Indeed, for every compactly supported smooth function 
$H \colon M \times [0,1] \to \RR$ each time-$t$ map~$\phi_H^t$ of its Hamiltonian flow
is a symplectomorphism.
Symplectomorphisms of this form are called Hamiltonian diffeomorphisms.
The Hamiltonian flow is the flow generated by the vector field~$X_H$ 
implicitely defined by
$$
\omega (X_H, \cdot) \,=\, -dH (\cdot) .
$$
For $(\RR^{2n}, \omega_0)$ one has $X_H = J_0 \2 \nabla H$, where $J_0$ is the usual
complex structure on $\oplus_i \RR^2 (x_i,y_i)$.

\section{Symplectic embedding problems}

By Darboux's Theorem, symplectic manifolds have no local invariants 
beyond the dimension. 
But there are several ways to associate global numerical invariants 
to symplectic manifolds. 
One of them is by looking at embedding problems.
Take a compact subset~$K$ of~$(\RR^{2n}, \omega_0)$.
By a symplectic embedding $K \to M$
we mean the restriction to~$K$ of a smooth embedding $\gf \colon U \to M$
of an open neighbourhood of~$K$ that is symplectic, $\gf^* \omega = \omega_0$.
In this case we write $K \se (M,\omega)$.
For every $K$, the largest number $\lambda$ such that the dilate~$\lambda K$
symplectically embeds into $(M,\omega)$ is then a symplectic invariant of $(M,\omega)$.
Seven further reasons to study symplectic embedding problems can be found 
in~\cite{Sch18}.

\subsection{The Nonsqueezing Theorem}  \label{ss:nonsq}
Now take the closed ball $\B^{2n}(a)$ of radius $\sqrt{a/\pi}$ centred at the origin of~$\RR^{2n}$.
(The notation reflects that symplectic measurements are 2-dimensional.)
By what we said above, there are very many symplectic embeddings
$\B^{2n}(a) \se \RR^{2n}$.
However, none of them can make the ball thinner, 
as Gromov proved in his pioneering paper~\cite{Gr85}:

\m \ni
{\bf Nonsqueezing Theorem.}
{\it $\B^{2n}(a) \se \B^2(A) \times \RR^{2n-2}$  \1 only if \1 $a \leq A$.}

\m
The identity embedding thus already provides the largest ball that 
symplectically fits into the cylinder $\B^2(A) \times \RR^{2n-2}$ of infinite volume. 
While there are many forms of symplectic rigidity, this theorem is its most fundamental manifestation.
The theorem shows that some volume preserving mappings cannot be approximated 
by symplectic mappings in the $C^0$-topology. 
It also gave rise to the useful concept of a symplectic capacity,~\cite{EkHo89}.

\s
In \cite{LaMc95}, Lalonde and McDuff generalized the Nonsqueezing Theorem to all symplectic manifolds:

\m \ni
{\bf General Nonsqueezing Theorem.}
{\it
For any symplectic manifold $(M,\omega)$ of dimension~$2n-2$, 
$$
\mbox{$\B^{2n}(a) \se \bigl( \B^2(A) \times M, \omega_0 \oplus \omega \bigr)$  \1 only if \1 $a \leq A$.}
$$
}

Every good tool in symplectic geometry can be used to prove the Nonsqueezing Theorem. 
However, the technique of $J$-holomorphic curves used by Gromov is the most influential
one, and also the most important tool in McDuff's work.

An almost complex structure~$J$ on a manifold~$P$ is a smooth collection $\{ J_p \}_{p \in P}$, 
where $J_p$ is a linear endomorphism of~$T_pP$ such that $J_p^2 = -\id$.
The `almost' indicates that such a structure does not need to be a complex structure, 
i.e.\ does not need to come from a holomorphic atlas. 
Not all symplectic manifolds admit complex structures, but they all admit almost complex structures.
A $J$-holomorphic curve in an almost complex manifold~$(P,J)$ is a map~$u$ 
from a Riemann surface~$(\Sigma, j)$ to~$(P,J)$ such that
$$
du \circ j \,=\, J \circ du.
$$
This equation generalizes the Cauchy--Riemann equation defining holomorphic maps $\CC \to \CC^n$.
In this text, the domain of a $J$-holomorphic curve will always be the usual Riemann sphere, 
namely the round sphere~$S^2 \subset \RR^3$ whose complex structure~$j$ rotates a vector $v \in T_pS^2$
by~$\frac \pi 2$. 
Even in this case, it is usually impossible to write down a $J$-holomorphic curve 
for a given~$J$.
But this is not a problem, since one usually just wants to know that such a curve exists.

Now assume that $\gf \colon \B^{2n}(a) \se \B^2(A) \times \RR^{2n-2}$.
Choose $k$ so large that after a translation the image of~$\gf$ is contained in 
$\B^2(A) \times (0,k)^{2n-2}$.
Compactifying the disc to the sphere~$S^2(A')$ with its usual area form of area~$A'>A$
and taking the quotient to the torus $T^{2n-2} = \RR^{2n-2} / k \ZZ^{2n-2}$, 
we then obtain a symplectic embedding
$$
\Phi \colon \B^{2n}(a) \se S^2(A') \times T^{2n-2} =: (P,\omega) 
$$
where $\omega$ is the split symplectic structure on the product~$P$.
We will see that $a \leq A'$. Since $A'>A$ was arbitrary, Gromov's theorem then follows.

Denote by $J_0$ the usual complex structure on $\B^{2n}(a) \subset \CC^n$,
and let $J$ be an almost complex structure on~$P$ that restricts
to $\Phi_* J_0$ on $\Phi (\B^{2n}(a))$ and that is {\it $\omega$-tame}, 
meaning that $\omega (v, Jv)>0$ for all nonzero $v \in TP$.
Such an extension exists, since 
$\omega$-tame almost complex structures can be viewed as sections
of a bundle over~$P$ whose fibers are contractible.

\m \ni
{\bf Lemma.} 
{\it 
There exists a $J$-holomorphic sphere~$u(S^2)$ through~$\Phi (0)$
that represents the homology class of~$S^2(A')$.
}

\m
This existence result follows from Gromov's compactness theorem for $J$-holomorphic curves 
in symplectic manifolds.
The key point for the proof of the compactness theorem is that $J$ is $\omega$-tame, 
implying that $J$-holomorphic curves cannot behave too wildly.
The compactness theorem implies the lemma because the class of~$S^2(A')$ is primitive in~$H_2(P;\ZZ)$.

\begin{figure}[h]
 \begin{center}
  \psfrag{gf}{$\Phi$}
  \psfrag{0}{$0$}
  \psfrag{S}{$S$}
  \psfrag{B}{$\textcolor{blue}{\B^{2n}(a)}$}
  \psfrag{C}{$T^{2n-2}$}
  \psfrag{u}{$\textcolor{Green1}{u(S^2)}$}
  \psfrag{f0}{$\Phi(0)$}
  \psfrag{fs}{$\textcolor{red}{\Phi(S)}$}
  \psfrag{fb}{$\textcolor{blue}{\Phi(\B^{2n}(a))}$}
  \psfrag{da}{$S^2(A')$}
  \leavevmode\epsfbox{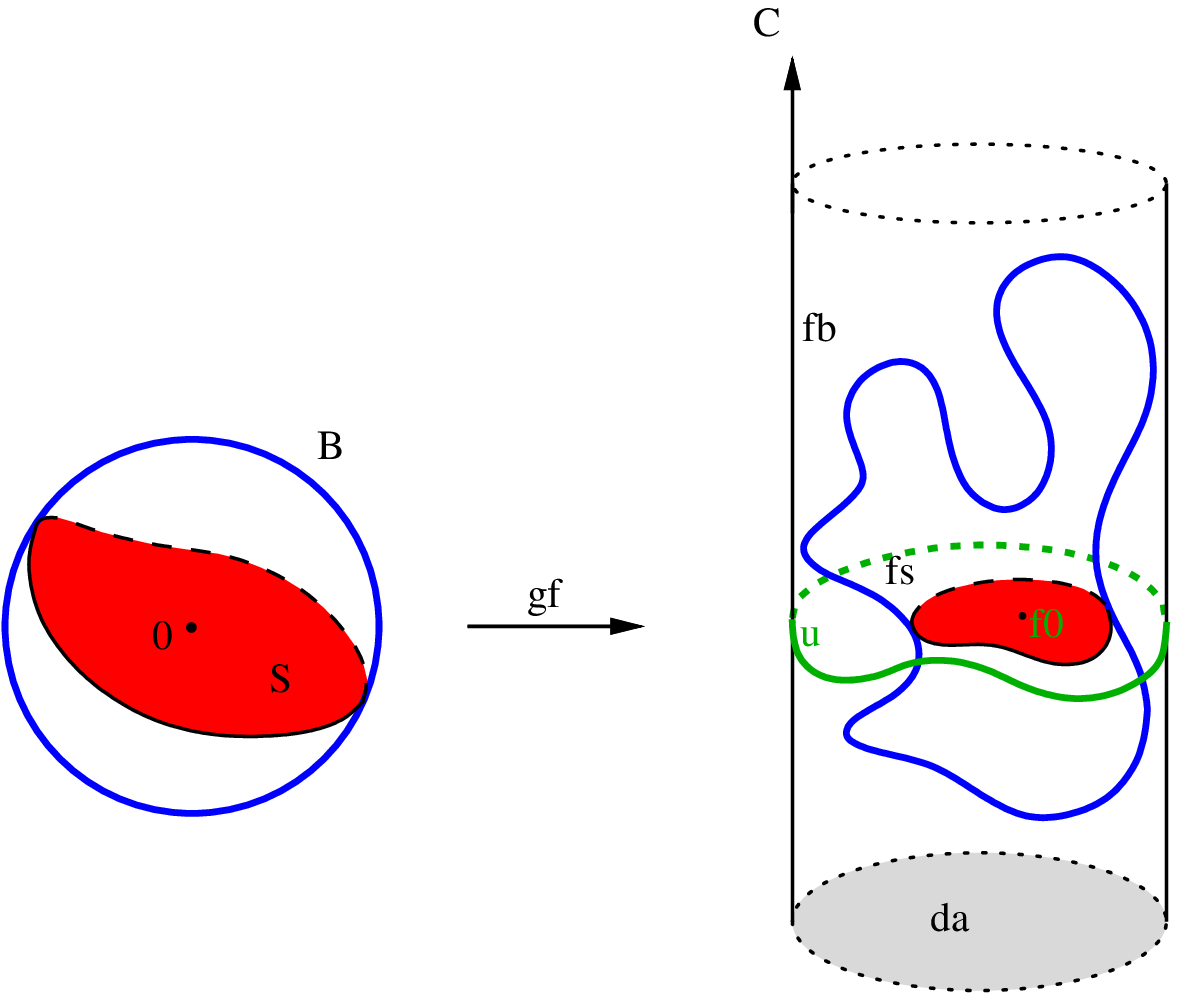}
 \end{center}
 \label{fig.nonsquez}
\end{figure}

The Nonsqueezing Theorem readily follows from the lemma:
The set
$$
\textcolor{red}{S} \,=\, \Phi^{-1} \left( \textcolor{red}{\Phi (\B^{2n}(a)) \cap u(S^2)} \right)
$$
is a proper 2-dimensional complex surface in~$\B^{2n}(a)$ through~$0$.
By the Lelong inequality from complex analysis, the area of~$S$ with respect to the
Euclidean inner product is at least~$a$.
Using also that $\Phi$ is symplectic we can now estimate
$$
a \,\leq\, \area (S) \,=\, \int_S \omega_0 \,=\, \int_S \Phi^*\omega \,=\, \int_{\Phi (S)} \omega
\,\leq\, \int_{u(S^2)} \omega \,=\, A'.
$$

\m
This proof of the Nonsqueezing Theorem also works, for instance, for closed symplectic manifolds~$(M,\omega)$ 
for which the integral of~$\omega$ over spheres vanishes. 
In the general case, however, there are troubles with holomorphic spheres of negative Chern number.
By now, these troubles are overcome thanks to the definition of Gromov--Witten invariants for 
general closed symplectic manifolds, 
by work of several (teams of) authors. 
McDuff helped to clarify the approach of 
Fukaya--Oh--Ono--Ohta in her joint work with Wehrheim, see~\cite{McWe17, 19}. 
In \cite{LaMc95}, however, Lalonde and McDuff circumvented 
all technical issues by using a chain of beautiful geometric constructions 
that reduced the General Nonsqueezing Theorem essentially to the case 
proved by Gromov.
The most important and influential of these is a multiple folding construction.
Its simple version
was used in~\cite{LaMc95} to show how the General Nonsqueezing Theorem implies 
that for all symplectic manifolds Hofer's metric on the group of compactly supported Hamiltonian diffeomorphisms 
is non-degenerate and hence indeed a metric, see~\S \ref{ss:Hofer}.
These two theorems were the first deep results in symplectic geometry proven for 
all symplectic manifolds.

\subsection{Ball packings} \label{ss:ball} 

We next try to pack a symplectic manifold by balls as densely as possible. 
Taking the open ball $\ballvier (1)$ as target, let
$$
p_k \,=\, \sup 
 \left\{ \frac{k \2 \Vol (\B^4(a))}{\Vol (\B^4(1))} \;\bigg|\, \coprod_k \B^4(a) \se 
\ballvier (1) \right\} 
$$
be the percentage of the volume of $\ballvier (1)$ that can be filled by $k$ symplectically embedded
equal balls. 
Then
\begin{equation} \label{t:2}
 \renewcommand{\arraystretch}{1.5}
  \begin{array}{c|ccccccccc} \hline
  k    & 1 &  2  &  3  & 4 & 5  & 6 & 7 & 8 & \geqslant 9      \\ \hline \hline
  p_k  & 1 & \frac 12 & \frac{3}{4} & 1  &  \frac{20}{25}  &  \frac{24}{25}  & \frac{63}{64} & \frac{288}{289}  & 1  \\ \hline
	c_k  & 1 & 2 & 2 &  2 & \frac 52  & \frac{5}{2} & \frac{8}{3} & \frac{17}{6} &
  \sqrt{k}  \\ \hline
  \end{array}
\end{equation}

\m \ni
The lower line gives the capacities
$$ 
c_k \,=\, \inf\, \biggl\{ A \;\bigg|\, \coprod_k \B^4(1) \se \ballvier (A) \biggr \} 
$$
that are related to the packing numbers $p_k$ by $c_k^2 = \frac{k}{p_k}$.
This table was obtained for $k \leqslant 5$ by Gromov~\cite{Gr85}, 
for $k=6,7,8$ and $k$ a square by McDuff and Polterovich~\cite{McPo94},
and for all~$k$ by Biran~\cite{Bi97}.
This result is a special case of the following 
algebraic reformulation of the general ball packing problem
\begin{equation} \label{e:Bai}
\coprod_{i=1}^k \B^4(a_i) \se \ballvier (A) .
\end{equation}

\m \ni
{\bf Ball Packing Theorem.}
{\it 
An embedding~\eqref{e:Bai} exists if and only if

\s
\begin{itemize}
\item[(i)]
{\rm (Volume constraint)}\, $A^2 > \sum_{i=1}^k a_i^2$;

\s
\item[(ii)]
{\rm (Constraint from exceptional spheres)}\,
$A > \frac 1d \sum_{i=1}^k a_i \1 m_i$ for every vector of non-negative integers $(d;m_1, \dots, m_k)$
that solves the Diophantine system
\begin{equation} \label{eq:ee}
\sum_i m_i = 3d-1, \qquad \sum_i m_i^2 = d^2+1 \tag{DE}
\end{equation}
and can be reduced to $(0;-1,0,\dots,0)$ 
by repeated Cremona moves.
\end{itemize}
}

\m \ni
Here, a Cremona move takes a vector $(d;m_1, \dots, m_k)$ with $m_1 \geqslant \dots \geqslant m_k$
to the vector 
$$
(d';\mm') = (d + \delta; m_1+\delta, m_2+\delta, m_3+\delta, m_4, \dots, m_k),
$$
where $\delta = d - (m_1+m_2+m_3)$, and then reorders $\mm'$.

\m
Before discussing the proof, we use the theorem to obtain Table~\ref{t:2}.
If $(d;m_1, \dots, m_k)$ is a solution of \eqref{eq:ee}, then 
$$
(3d-1)^2 \,=\, \left( \sum_{i=1}^k m_i \right)^2 \,\leq\, k \sum_{i=1}^k m_i^2 \,=\, k \2 (d^2+1) ,
$$
that is, 
$$
(9-k) \2 d^2 - 6d + (1-k) \,\leq\, 0 .
$$
For $k \leq 8$, this equation has finitely many solutions $d$,
and so \eqref{eq:ee} has finitely many solutions for $k \leq 8$.
They are readily computed:
$$
(1; 1, 1),\;\; (2; 1^{\times 5}),\;\;
(3; 2, 1^{\times 6}), \;\; 
(4;2^{\times3}, 1^{\times5}),\;\;
(5; 2^{\times6}, 1,1),\;\; (6; 3, 2^{\times7}),
$$
and all these vectors reduce to $(0;-1)$ by Cremona moves.
For instance, for the problem $\coprod_8 \B^4(1) \se \ballvier (A)$ 
the strongest constraint comes from the solution 
$(6; 3, 2^{\times7})$, that gives $A > \frac{17}{6}$.

On the other hand, if $(d;m_1, \dots, m_k)$ is a solution of~\eqref{eq:ee} with $k \geq 9$, then
$$
\frac ad \2 \sum_{i=1}^k m_i \,=\, \frac ad \, (3d-1) \,<\, \frac ad \, 3d \,\leq\, 
\frac ad \, \sqrt k \, d \,=\, a \, \sqrt k .
$$
Hence the constraint $A > \frac ad \sum_{i=1}^k m_i$ is weaker than 
the volume constraint $A^2 > k \1 a^2$, and so $p_k=1$.

\m 
The proof of the Ball Packing Theorem is a beautiful story in three chapters, 
each of which contains an important idea of McDuff.
The original symplectic embedding problem is converted 
to an increasingly algebraic problem in three steps. 

The starting point is the symplectic blow-up construction, 
that goes back to Gromov and Guillemin--Sternberg, and was first used by 
McDuff~\cite{91a} to study symplectic embeddings of balls.
Recall that the complex blow-up $\Bl (\CC^2)$ of $\CC^2$ at the origin~$0$
is obtained by replacing~$0$
by all complex lines in~$\CC^2$ through~$0$.
At the topological level, this operation can be done as follows:
First remove from $\CC^2$ an open ball $\ballvier$. 
The boundary $S^3$ of $\CC^2 \setminus \ballvier$ is foliated by the Hopf circles
$\{ (\alpha z_1, \alpha z_2) \mid \alpha \in S^1 \}$, namely the 
intersections of $S^3$ with complex lines.
Now $\Bl (\CC^2)$ is obtained by replacing each such circle by a point. 
The boundary sphere~$S^3$ becomes a 2-sphere~$\CP^1$ in~$\Bl (\CC^2)$
of self intersection number~$-1$, called the exceptional divisor.
The manifold~$\Bl (\CC^2)$ is diffeomorphic to the connected sum~$\CC^2 \# \overline \CP^2$.

This construction can be done in the symplectic setting:
If one removes $\ballvier (a)$ from~$\RR^4$, 
then there exists a symplectic form $\omega_a$ on~$\Bl (\RR^4)$ such that $\omega_a = \omega_0$ 
outside a tubular neighbourhood of the exceptional divisor~$\CP^1$ and such that $\omega_a$ is symplectic
on~$\CP^1$ with $\int_{\CP^1} \omega_a = a$. 
Given a symplectic embedding $\gf \colon \B^4(a) \to (M,\omega)$ into a symplectic 4-manifold,
we can apply the same construction to $\gf (\B^4(a))$
in $M$ to obtain the symplectic blow-up of $(M,\omega)$ by weight~$a$.

There is also an inverse construction:
Given a symplectically embedded $-1$~sphere~$\Sigma$ in a symplectic $4$-manifold~$(M,\omega)$ 
of area $\int_\Sigma \omega = a$, one can cut out a tubular neighbourhood of~$\Sigma$
and glue back~$\B^4(a)$, to obtain the `symplectic blow-down' of~$M$.

Now let $M_k$ be the smooth manifold obtained by blowing up the complex projective plane~$\CP^2$ in $k$~points.
Its homology $H_2(M_k;\ZZ)$ is generated by the class $L$ of a complex line and by the classes 
$E_1, \dots, E_k$ of the exceptional divisors.
Let $\ell, e_i \in H^2(M_k;\ZZ)$ be their Poincar\'e duals.
Every symplectic form $\omega$ on~$M_k$ defines a first Chern class~$c_1(\omega)$,
namely the first Chern class of any $\omega$-tame almost complex structure.
If we take a symplectic form~$\omega$ on~$M_k$ constructed as above via $k$ symplectic ball embeddings
(that always exist if the balls are small enough), 
then $-c_1(\omega)$ is Poincar\'e dual to the class $K := -3L+\sum_i E_i$.
Define $\cc_K (M_k) \subset H^2(M_k;\RR)$ to be the set of classes represented by symplectic forms
with $-c_1(\omega) = \PD (K)$.

Compactifying $\ballvier (A)$ to~$\CP^2 (A)$ 
with its usual K\"ahler form integrating to~$A$ over a complex line, 
and using that the classes~$E_i$ can be represented by symplectic $-1$~spheres, 
McDuff and Polterovich~\cite{McPo94} obtained

\m \ni
{\bf Step 1.}
{\it
There exists an embedding $\coprod_{i=1}^k \B^4(a_i) \se \ballvier (A)$ 
if and only if 
the class 
$\ga := A \ell - \sum_{i=1}^k a_i \1 e_i \in H^2(M_k;\RR)$ 
lies in $\cc_K (M_k)$.
}

\m
Our symplectic embedding problem is thus translated into a problem on the symplectic cone~$\cc_K (M_k)$.
In K\"ahler geometry, the problem of deciding 
which cohomology classes can be represented by a K\"ahler form 
has a long history and is quite well understood.
The solution of our symplectic analogue, however, needs different ideas and tools:
Call a class $E \in H_2(M_k;\ZZ)$ exceptional if $K \cdot E = 1$ and $E^2 = -1$, 
and if $E$ can be represented by a smoothly embedded $-1$~sphere.

If $\omega$ is a symplectic form on $M_k$ with $-c_1(\omega) = \PD (K)$,
then any $\omega$-symplectic embedded $-1$~sphere represents an exceptional class.
Seiberg--Witten--Taubes theory implies that the converse is also true:
Every exceptional class~$E$ can be represented by an $\omega$-symplectic embedded $-1$~sphere, \cite{LiLi}.
This implies one direction of 

\m \ni
{\bf Step 2.}
{\it
$\ga = A \ell - \sum_i a_i \1 e_i \in H^2(M_k;\RR)$ 
lies in $\cc_K (M_k)$
if and only if $\ga^2 >0$ and $\ga (E) >0$ for all exceptional classes.
}

\m
This equivalence is remarkable: 
Of course, a necessary condition for a class~$\ga$ with positive square
to have a symplectic representative is that $\ga$ evaluates positively on all classes that can be represented by closed symplectically embedded surfaces (and in particular on spheres). 
But the equivalence says that this is also a sufficient condition, 
and that it is actually enough to check positivity on spheres.

To prove the other direction in Step~2, 
one starts with an embedding of $k$ tiny balls of size $\gve \1 a_i$
and then changes the symplectic form on $M_k$ in class $A \ell - \gve \sum_i a_i \1 e_i$
in such a way that these balls look large.
This can be done with the help of the inflation method of Lalonde--McDuff~\cite{94, LaMc96a}.

\m
\ni 
{\bf Inflation Lemma.}
{\it
Let $(M,\go)$ be a closed symplectic 4-manifold, 
and assume that the class $C \in H_2(M;\ZZ)$ with $C^2 \geqslant 0$
can be represented by a closed connected embedded $J$-holomorphic curve~$\Sigma$ for some $\go$-tame~$J$.
Then the class $[\go] + s \PD (C)$ has a symplectic representative 
for all $s \geqslant 0$.
}

\begin{figure}[ht] 
 \begin{center}
  \psfrag{P}{$\PD(C)$}
  \psfrag{om}{$[\omega]$}
  \psfrag{H}{$H^2(M;\RR)$}
  \leavevmode\epsfbox{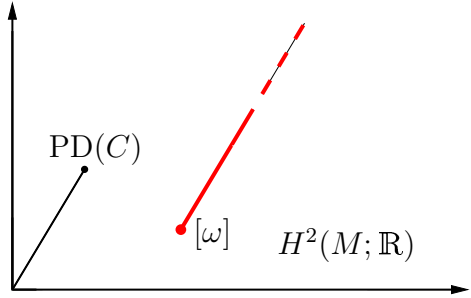}
 \end{center}
 \caption{The ray in the symplectic cone provided by the Inflation Lemma} 
 \label{fig:inflationlemma}
\end{figure}
%

\m
I explain the proof for the case $C^2 = 0$.
In this case, the normal bundle of~$\Sigma$ is trivial, and we can
identify a tubular neighbourhood of~$\Sigma$ with $\Sigma \times D$,
where $D \subset \RR^2$ is a disc.
Pick a radial function $f(r)$ with support in~$D$ that is non-negative and has $\int_D f = 1$.
Let $\beta$ be the closed 2-form on~$M$ that equals
$\beta (z,x,y) = f(r) \, dx \wedge dy$ on $\Sigma \times D$ and vanishes outside of $\Sigma \times D$.
Then $[\beta] = \PD ([\Sigma])$ and the forms 
$\omega + s \beta$ are symplectic for all $s \geq 0$. Indeed, 
$$
(\omega + s \beta)^2 \,=\, \underbrace{\omega^2}_{>0} \,+\, 
     2 s \underbrace{\omega \wedge \beta}_{\geq 0} \,+\, s^2\underbrace{\beta^2}_{= 0} \,>\, 0
$$
where for the middle term we used that $\omega |_\Sigma$ is symplectic.

Now take an embedding $\coprod_{i=1}^k \B^4(\gve \1 a_i) \se \ballvier (A)$
of tiny balls.
By Step~1 we know that the class $\alpha_{\gve} := A \ell - \gve \sum_i a_i \1 e_i$
has a symplectic representative~$\omega_{\gve}$.
We wish to inflate this form to a symplectic form in class~$\alpha$.
A first try could be to inflate $\omega_{\gve}$ directly in the direction $\alpha - \alpha_{\gve}$
to get up to~$\alpha$.
But this does not work, because
$$
(\alpha-\alpha_{\gve})^2 \,=\, \bigl( - \textstyle{\sum}_i (1-\gve) \2 a_i \2 e_i \bigr)^2 \,=\, 
-(1-\gve)^2 \, \textstyle{\sum}_i \2 a_i^2 \,<\, 0.
$$
However, assuming for simplicity that $A$ and the~$a_i$ are rational, 
Seiberg--Witten--Taubes theory implies that there exists an integer~$n$ such
that the Poincar\'e dual of $n \alpha \in H^2 (M_k;\ZZ)$ can be represented
by a connected embedded $J$-holomorphic curve for a generic $\omega$-tame~$J$.
We can thus inflate $\omega$ in the direction of~$n \alpha$ and obtain that
the classes $\alpha_\gve + s \1 n  \1 \alpha$ have symplectic representatives for all $s \geq 0$.
Rescaling these forms by $\frac{1}{sn+1}$, we obtain symplectic forms in classes 
$A\ell - \frac{sn+\gve}{sn+1} \sum_i a_i \1 e_i$ that are as close to~$\alpha$ as we like.

\begin{figure}[ht] 
 \begin{center}
  \psfrag{1}{$1$}
  \psfrag{A}{$A$}
  \psfrag{l}{$\ell$}
  \psfrag{a}{$\ga$}
  \psfrag{ae}{$\ga_\gve$}
  \psfrag{sum}{$-\sum a_i\,e_i$}
  \psfrag{e}{$\gve$}
  \psfrag{na}{$\ga_{\gve} + s n\ga$}
  \psfrag{N}{\mbox{No!}}
  \leavevmode\epsfbox{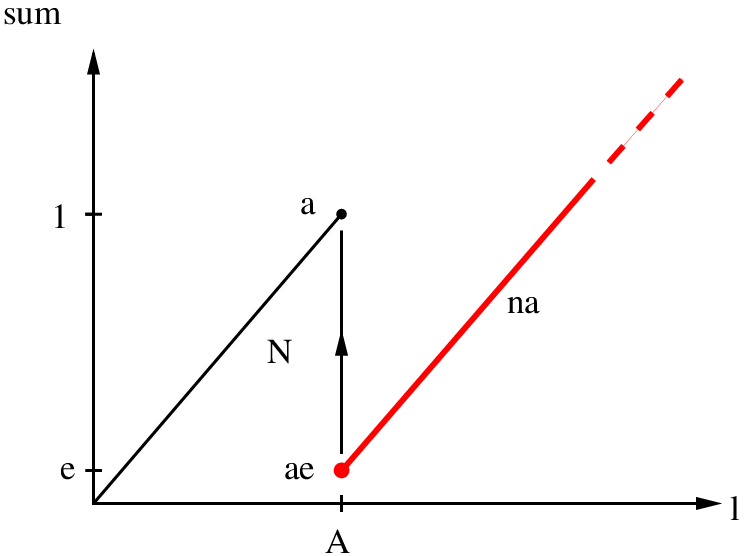}
 \end{center}
 \caption{Inflating $\ga_\eps = [\go_\eps]$ to~$\ga$} 
 \label{fig:inflation}
\end{figure}

\m
Usually, $J$-holomorphic curves are used to {\it obstruct}\/ certain symplectic embeddings, 
like in the Nonsqueezing Theorem. 
In the Inflation Lemma and Step~2, however, other $J$-holomorphic curves are
used to {\it construct}\/ symplectic embeddings.

\m \ni
{\bf Step 3.}
{\it
$\ga^2 >0$ and $\ga (E) >0$ for all exceptional classes
if and only if (i) and (ii) in the theorem hold.
}

\m
Much of this step is rewriting:
The inequality $\alpha^2 >0$ translates to the volume constraint~(i),
and the the inequality $\ga (E) >0$ for an exceptional class 
translates to~(ii) if we write~$E$ in the natural basis, 
$E = d L - \sum_{i=1}^k m_i \1 E_i$.
We are thus left with showing that a class $E = (d;m_1,\dots,m_k)$ 
that satisfies~\eqref{eq:ee}
is exceptional exactly if it reduces to~$(0;-1)$ under Cremona moves. 
This follows from a combinatorial argument~in \cite{LiLi}, see~\cite{McSch12}.

\subsection{Ellipsoids}

We now look at embeddings of 4-dimensional ellipsoids
$$
\E (a,b) \,=\,
\biggl\{ (z_1, z_2) \in \CC^2 \mid \frac{\pi |z_1|^2}{a} + \frac{\pi |z_2|^2}{b} \leq 1 \biggr\} .
$$
This is the ellipsoid in~$\CC^2$ whose projections to the coordinate planes are closed
discs of area~$a$ and~$b$.
Again we take as target a ball, and after rescaling study the function
$$
c(a) \,=\, \inf \left\{ A \mid \E (1,a) \se \ballvier (A) \right \}, \qquad a \geqslant 1.
$$
Since this function is continuous, we can assume that $a$ is rational. 
Then each leaf of the characteristic foliation on the boundary of $\E(1,a)$ is closed.
Proceeding as before, we compactify $\ballvier (A)$
to~$\CP^2(A)$ and, given an embedding $\E(1,a) \se \ballvier (A) \subset \CP^2(A)$,
remove the image and collapse the remaining boundary along the characteristic foliation.
But now this yields a symplectic orbifold with one or two cyclic quotient singularities, 
coming from the special leaves in the two coordinate planes.
It may be difficult to reprove the results from 
Seiberg--Witten--Taubes theory used in the last section in such a space.
McDuff in~\cite{09} simply circumvented the singularities by
using a version of the Hirzebruch--Jung resolution of singularities.
She removed a bit more than the ellipsoid by successively blowing up finitely many balls,
thereby producing a chain of $J$-spheres. 
Inflating this chain she reduced the problem $\E(1,a) \se \ballvier  (A)$ 
to the ball packing problem~\eqref{e:Bai}:
\begin{equation} \label{e:equiv}
\E(1,a) \se  \ballvier  (A)
\;\Longleftrightarrow\;\, 
\coprod_{i=1}^k \B^4(a_i) \se \ballvier  (A)
\end{equation}
where the $a_i$ are given by
$$
(a_1, \dots, a_k) \,=:\, \ww (a) \,=\, \bigl(\underbrace{1,\dots,1}_{\ell_0}, \,
\underbrace{w_1,\dots,w_1}_{\ell_1}, \,
\dots, \, \underbrace{w_N,\dots,w_N}_{\ell_N} \bigr)
$$
with the weights $w_i >0$
such that $w_1 = a-\ell_0 < 1$, $w_2 = 1-\ell_1 w_1 < w_1$, and so on.
For instance, $\ww (3) = (1,1,1)$ and
$\ww \bigl(\tfrac{11}{4}\bigr) = \left( 1,1, \tfrac 34, \tfrac 14, \tfrac 14, \tfrac 14 \right)$.
The multiplicities~$\ell_i$ of~$\ww (a)$ give the continued fraction expansion of~$a$.
For $a \in \NN$, \eqref{e:equiv} specializes to
\begin{equation} \label{e:equivk}
\E(1,k) \se \ballvier  (A) \;\Longleftrightarrow\;\, 
\coprod_k \B^4(1) \se 
\ballvier  (A) .
\end{equation}
The ball packing problem $\coprod_k \B^4(1) \to \ballvier  (A)$
is thus included in the 1-parameter problem $\E(1,a) \se \ballvier (A)$.

\begin{figure}[ht]
 \begin{center}
  \psfrag{a}{$a$}
  \psfrag{wa}{$\sqrt{a}$}
  \psfrag{an}{$a_n$}
  \psfrag{an1}{$a_{n+1}$}
  \psfrag{bn}{$b_n$}
  \psfrag{c}{$\textcolor{red}{c(a)}$}
  \psfrag{gn}{$\gg_n$}
  \psfrag{gn1}{$\gg_{n+1}$}
 \leavevmode\epsfbox{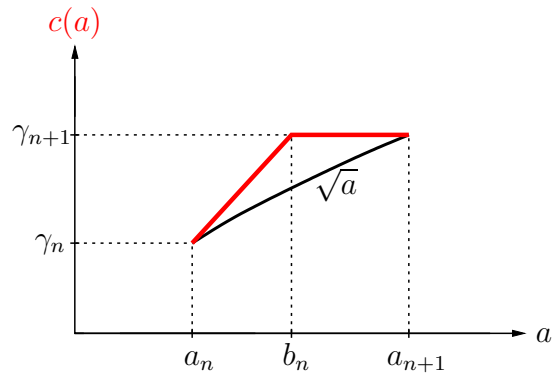}
 \end{center}
 \caption{The $n^{\mbox{\scriptsize th}}$ step of the Fibonacci stairs}
 \label{fig.fib.step}
\end{figure}

The function $c(a)$ was computed in~\cite{McSch12} with the help of~\eqref{e:equiv}.
The volume constraint is now $c(a) \geqslant \sqrt{a}$.
Recall that the Fibonacci numbers are recursively defined by 
$f_{-1}=1, f_0=0, f_{n+1} = f_n + f_{n-1}$.
Denote by $g_n: = f_{2n-1}$ the odd-index Fibonacci numbers,
hence $(g_0, g_1, g_2, g_3, g_4,\dots ) = (1,1,2,5,13,\dots)$.
The sequence $\gg_n := \frac{g_{n+1}}{g_n}$, 
$$
\left( \gg_0, \gg_1, \gg_2, \gg_3, \dots \right) \,=\,
\left( 1, 2, \frac 52, \frac{13}{5}, \dots \right) ,
$$
converges to $\tau^2$, where $\tau := \frac{1+\sqrt{5}}{2}$ is the Golden Ratio.
Define the {\it Fibonacci stairs}\/ as the graph on~$[1,\tau^4]$ 
made from the infinitely many steps shown in Figure~\ref{fig.fib.step},
where $a_n = \gg_n^2 = \bigl( \frac{g_{n+1}}{g_n} \bigr)^2$ and $b_n = \frac{g_{n+2}}{g_n}$.
The slanted edge starts on the volume constraint $\sqrt{a}$ and extends to a line through the origin.

\begin{figure}[ht]
 \begin{center}
  \psfrag{1}{$1$}
  \psfrag{2}{$2$}
  \psfrag{4}{$2^2$}
  \psfrag{5}{$5$}
  \psfrag{t}{$(\tau^4,\tau^2)$}
  \psfrag{a}{$a$}
  \psfrag{c}{$\textcolor{red}{c(a)}$}
  \psfrag{52}{$\tfrac 52$}
  \psfrag{13}{$\tfrac{13}{5}$}
  \psfrag{25}{$(\tfrac{5}{2})^2$}
  \psfrag{132}{$(\tfrac{13}{5})^2$}
 \leavevmode\epsfbox{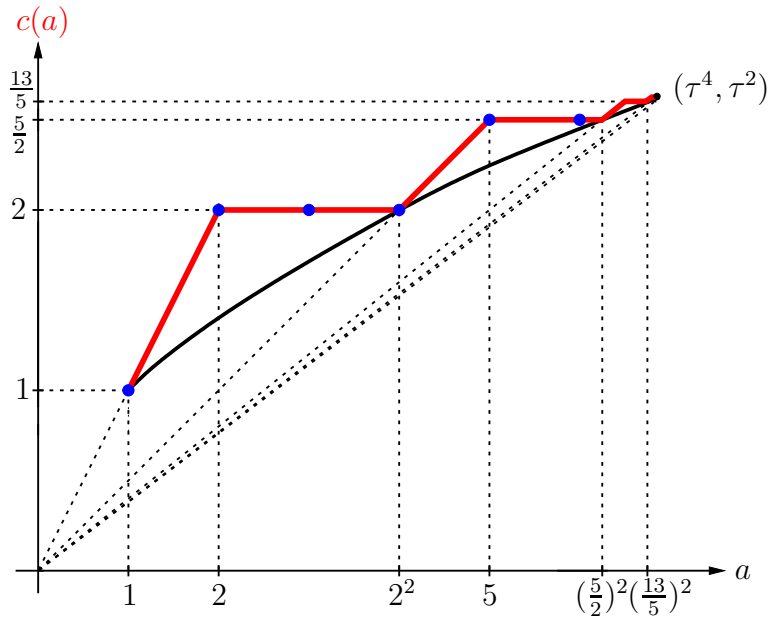}
 \end{center}
 \caption{The Fibonacci stairs: The graph of $c(a)$ on $\left[ 1,\tau^4 \right]$}
 \label{figure.Fib}
\end{figure}
%
%

%
%
%

\m
\ni
{\bf Ellipsoid Embedding Theorem.} 
{\it

\s
\begin{itemize}
\item[(i)]
On the interval $\left[1,\tau^4\right]$ the function $c(a)$ is given 
by the Fibonacci stairs. 

\m
\item[(ii)]
On the interval $\left[\tau^4, (\frac{17}{6})^2 \right]$ we have
$c(a)=\sqrt a$ except on nine disjoint intervals where
$c$ is a step made from two segments. 
The first of these steps has the vertex at~$(7, \frac 83)$ and the last 
at~$(8,\frac{17}{6})$. 

\m
\item[(iii)]
$c(a)=\sqrt a$\, for all $a \geqslant (\frac{17}{6})^2$.
\end{itemize}
}

\m
Thus the graph of $c(a)$ starts with an infinite completely regular staircase, then
has a few more steps, but for $a \geqslant (\frac{17}{6})^2 = 8 \frac{1}{36}$
is given by the volume constraint.

\begin{figure}[h] 
 \begin{center}
  \psfrag{a}{$a$}
  \psfrag{1}{$1$}
  \psfrag{tt}{$\tau^4$}
  \psfrag{8}{$(\frac{17}{6})^2$} 
  \psfrag{s}{\textcolor{red}{\mbox{structured rigidity}}}
  \psfrag{t}{\textcolor{magenta}{\mbox{transition}}}
  \psfrag{f}{\textcolor{blue}{\mbox{flexibility}}}
 \leavevmode\epsfbox{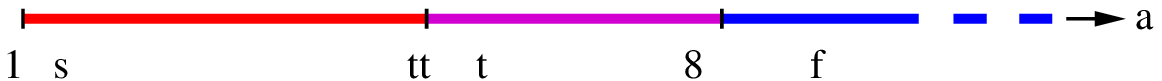}
 \end{center}
\end{figure}

The theorem better explains the packing numbers~$c_k$ in Table~\ref{t:2}
in view of the equivalence~\eqref{e:equivk}, cf.\ the blue dots in Figure~\ref{figure.Fib}.
The embedding constraint at~$b_n$ comes from
a really exceptional exceptional sphere, namely one in an exceptional class 
$E=(d;\mm)$ such that $\mm$ is parallel to the weight expansion~$\ww (b_n)$.

\m
When McDuff first looked at ellipsoid embeddings in~\cite{09},
it was not clear at all that this would unveil an interesting 
fine structure of symplectic rigidity. 
Many more infinite staircases have now been found, see for example \cite{FrMu15, Us19, CHMP20, BHMMMPW20}.
Usher's results in~\cite{Us19} show that 
a simple looking problem, 
such as that for which~$A$ the ellipsoid~$\E(1,a)$ embeds into 
the polydisc $\B^2(A) \times \B^2(bA)$, 
is already so intricate that we shall probably never know the answer for all~$b$.
McDuff in~\cite{11a} and Hutchings in~\cite{Hu11b}
proved that
\begin{equation} \label{e:Hofer}
\E(a,b) \se \E(c,d) \;\Longleftrightarrow\; N_k(a,b) \leqslant N_k(c,d) \;\mbox{ for all } k 
\end{equation}
where $(N_k(a,b))$ is the non-increasing sequence obtained by ordering the
set $\{ ma +nb \mid m,n \in \ZZ_{\geq 0} \}$.
Using his embedded contact homology (ECH),  
Hutchings~\cite{Hu11a} had associated with every starshaped  
subset $K \subset \RR^4$ a sequence of numbers~$c_k(K)$ that are
monotone with respect to symplectic embeddings, 
and for an ellipsoid these ECH~capacities equal
the above sequence, $c_k(\E (a,b)) = N_k(a,b)$.
The McDuff--Hutchings Theorem~\eqref{e:Hofer} therefore implied that 
ECH-capacities are a complete set of invariants 
for the problem of embedding one four-dimensional ellipsoid into another. 

\s
These results are all in dimension four, 
and until recently, not much was known about higher dimensional symplectic embedding problems
beyond the Nonsqueezing Theorem.
The reason is that in dimension~four, $J$-holomorphic curves are a much 
more powerful tool, 
because of positivity of intersections of such curves, and because one usually finds them
with the help of the 4-dimensional Seiberg--Witten--Taubes theory.
However, 4-dimensional ellipsoid embeddings also opened the door for understanding certain
symplectic embedding problems in higher dimensions. 
I describe how they led to packing stability in all dimensions.
Given a connected symplectic manifold~$(M,\omega)$ of finite volume,
let $p(M,\omega)$ be the smallest number (or infinity) such that 
for every $k \geq p(M,\omega)$ an arbitrarily large percentage of the volume of $(M,\omega)$ 
can be covered by $k$ equal symplectically embedded balls.
Table~\ref{t:2} shows that $p(\B^4) = 9$.
The finiteness of $p$ is now known for many symplectic manifolds, 
and in particular for balls in all dimensions. 
In this case, the ingredients of the proof are:

\s
\begin{itemize}
\item[(1)]
McDuff's observation from \cite{09} that the ellipsoid $\E (1,\dots,1,k)$ can be cut into $k$ equal balls.

\s
\item[(2)]
The Ellipsoid Embedding Theorem for $\E(1,a) \se \ballvier (A)$ 
and the McDuff--Hutchings Theorem~\eqref{e:Hofer}.

\s
\item[(3)]
Ellipsoids admit a suspension construction, \cite{BuHi11}: 
For any vectors $\boldsymbol{a} , \boldsymbol{b}, \boldsymbol{c}$,
$$
\E(\boldsymbol{a}) \se \E(\boldsymbol{b}) \;\Longrightarrow\; 
\E(\boldsymbol{a}, \boldsymbol{c}) \se \E(\boldsymbol{b}, \boldsymbol{c}) .
$$
\end{itemize}

\ni
For instance, (2) yields that 
$$
\E (1,k) \;\textcolor{magenta}{\se}\; \E (k^{1/3}, k^{2/3}) \quad\mbox{ and }\quad 
\E (1,k^{2/3}) \;\textcolor{blue}{\se}\; \B^4 (k^{1/3}) 
\qquad \mbox{for all $k \geqslant 21$,}
$$
see~\cite{BuHi13}.
Together with (1) and (3) we thus obtain that for these $k$
$$
\coprod_k \B^6(1) \,\se\, \E (1,1,k) \;\textcolor{magenta}{\se}\;\E (1,k^{1/3},k^{2/3}) \;\textcolor{blue}{\se}\; 
\B^6(k^{1/3}) 
$$
and hence $p(\B^6) \leq 21$.
Already Gromov~\cite{Gr85} proved that $p(\B^6) \geq 8$. Is it true that $p(\B^6) = 8$\,?

A symplectic manifold $(M,\omega)$ is {\it rational}\/ if
a multiple of $[\omega]$ takes rational values on all integral 2-cycles.
All closed rational symplectic manifolds have packing stability, \cite{BuHi13}.
The additional ingredient in the proof is that every such manifold 
can be completely filled by an ellipsoid.

\subsection{Connectivity} 

For a compact set $K \subset \RR^{2n}$ and a $2n$-dimensional 
symplectic manifold $(M,\omega)$
let $\Emb_\omega(K,M)$ be the space of symplectic embeddings
$K \to (M,\omega)$, with the $C^\infty$-topology.
The results discussed above tell us for some $K$ and~$(M,\omega)$
whether this space is empty or not. 
In the latter case one may study its topology. 
The first task is to see whether $\Emb_\omega(K,M)$ is connected or not.

\begin{figure}[h] 
 \begin{center}
  \leavevmode\epsfbox{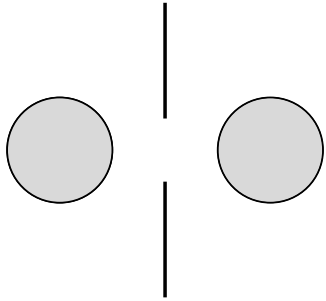}
 \end{center}
 \caption{Non-equivalent balls in the camel space} 
 \label{camel}
\end{figure}

We first take $K$ to be a ball $\B^{2n}(a)$.
Then $\Emb_\omega \bigl( \B^{2n}(a),M \bigr)$ need not be connected.
The first counter-example was Gromov's camel theorem in~\cite{Gr85}:
For $2n \geq 4$
the camel space in~$\RR^{2n}$ with eye of width~$1$ is the set 
$$
\cc^{2n}\,=\, 
\left\{ x_1<0 \right\} \cup \left\{ x_1>0 \right\} \cup \ballzwein (1).
$$ 
Now take $a>1$ and define the two embeddings $\gf_\pm \colon \B^{2n}(a) \to \cc^{2n}$
by $\gf_\pm(z) = z \pm v$, where $v$ is a multiple of $\frac{\pp}{\pp x_1}$
such that $\gf_+$ takes $\B^{2n}(a)$ to the right half space
and $\gf_-$ takes $\B^{2n}(a)$ to the left half space,
see Figure~\ref{camel}.
Then $\gf_+$ and $\gf_-$ are not isotopic, see~\cite{EG89, Vit92, MT93} for proofs.

On the other hand, McDuff showed in~\cite{91a} that 
$\Emb_\omega \bigl( \B^4(1), \ballvier (A) \bigr)$ is connected for all~$A>1$
and extended this result in~\cite{98, 09} to embeddings of any collection 
of 4-ellipsoids into a 4-ellipsoid.
For yet a further generalization see~\cite{Cr19}.

The key to these results is, again, symplectic inflation.
I give the idea of the proof for embeddings of one ball~$\B^4(1)$.
So assume we have two embeddings  
$$
\gf_1, \gf_2 \colon \B^4(1) \,\se\, \ballvier (A) .
$$
Since $\gf_1$ and~$\gf_2$ are close to linear maps near the origin, 
we can find $\gve >0$ and a compactly supported symplectic isotopy~$\psi$ 
of~$\ballvier (A)$ such that $\gf_1 = \psi \circ \gf_2$ on~$\B^4(\gve)$.
We may thus assume from the start that $\gf_1 = \gf_2$ on~$\B^4(\gve)$.

\begin{figure}[h] 
 \begin{center}
     \psfrag{0}{$0$}
     \psfrag{1}{$1$}
     \psfrag{2}{$2$}
     \psfrag{s}{$s$}
     \psfrag{e}{$\gve$}
     \psfrag{A}{$A \ell - f(s) \2 e$}
  \leavevmode\epsfbox{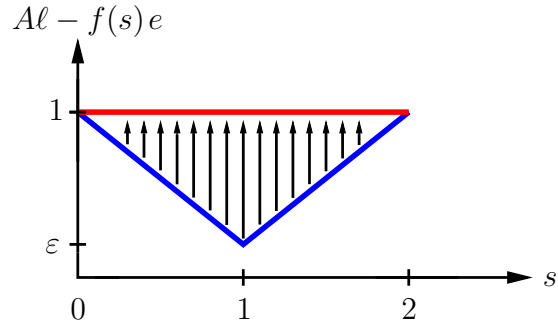}
 \end{center}
 \caption{Inflating $\omega_s$ to $\tilde \omega_s$}  
 \label{fig.inflation}
\end{figure}

Now consider the two paths of symplectic embeddings 
$\gf_j^t \colon \B^4(t) \to \ballvier (A)$, 
$\gve \leq t \leq 1$, defined by restricting the embeddings~$\gf_j$.
Reparametrizing the inverse of the path~$\gf_1^t$ on~$s \in [0,1]$ 
and the path~$\gf_2^t$ on~$s \in [1,2]$
we obtain a path of symplectic embeddings into~$\ballvier (A)$
that connects $\gf_1$ and~$\gf_2$.
After compactifying $\ballvier (A)$ to~$\CP^2(A)$ by adding the line~$\CP^1$
and by symplectically blowing up the images, this path of embeddings 
gives rise to a path of 
symplectic forms~$\omega_s$ on~$\Bl (\CP^2)$ in class~$\textcolor{blue}{[\omega_s]}$ 
as shown in Figure~\ref{fig.inflation}. 
As in~\S \ref{ss:ball} apply symplectic inflation to each form~$\omega_s$,  
deforming $\omega_s$ to a form~$\tilde \omega_s$ in class $[\omega_0] = [\omega_2] = \textcolor{red}{A \ell -e}$. 
Finally blow down the exceptional divisor for each~$s$ to get a path of embeddings 
$\phi_s \colon \B^4(1) \se \CP^2(A)$. 
Since each exceptional divisor is disjoint from~$\CP^1$,
each ball $\phi_s (\B^4(1))$ lies in $\CP^2(A) \setminus \CP^1$,
and by a theorem of McDuff from~\cite{90} this set with the symplectic form 
obtained after blow-down is indeed symplectomorphic to~$\ballvier (A)$,
cf.\ \S \ref{ss:structure} below.
Furthermore, $\phi_s$ connects $\gf_1$ with~$\gf_2$.

\s
McDuff's theorem that $\Emb_\omega \bigl( \B^4(1), \ballvier (A) \bigr)$
is connected is sharp in the following sense:
There are convex subsets $K$ of~$\RR^4$ with smooth boundary that are arbitrarily close to a ball 
and
such that $\Emb_{\omega}(K, \lambda \Kcirc)$ has at least two components, for some~$\lambda >1$, 
see~\cite{GuUs19}.
For the cube $\C^4(1) = \B^2(1) \times \B^2(1)$ the space
$\Emb_{\omega} \bigr(\C^4(1), \ballvier (3) \bigr)$ 
even has infinitely many connected components,~\cite{BMS20}.


\section{Other contributions to symplectic geometry}

\subsection{Construction of (counter-)examples}

McDuff again and again has given explicit constructions that are 
both elementary and profound.
You can grasp them quickly, yet they 
can be used to do many other things. 
We already encountered some of them in the previous section, 
and I now mention three more.

\subsection*{A simply-connected symplectic manifold that is not K\"ahler}

While K\"ahler manifolds are symplectic, the converse is not true.
For instance, the odd-index Betti numbers of compact K\"ahler manifolds are even. 
The first example of a compact symplectic manifold that is not K\"ahler
was found by Thurston; his example $(M,\omega)$ is a $T^2$-bundle over $T^2$ with $b_1(M) = 3$.
In \cite{84}, McDuff constructed simply-connected examples.
She symplectically embedded the Thurston manifold $M$ into $\CP^5$
and showed that the symplectic blow-up of~$\CP^5$ along~$M$ 
is simply-connected and has $b_3 = b_1(M) = 3$.

\subsection*{Cohomologous symplectic forms that are not diffeomorphic}

Consider two symplectic forms $\omega_0$ and $\omega_1$ on a closed symplectic manifold
that are cohomologous.
If these two forms can be connected by a smooth path $\omega_t$ of cohomologous symplectic forms,
then there exists an isotopy of~$M$ deforming $\omega_0$ to~$\omega_1$, by Moser's trick.
This is the case, for instance, for any two area forms on a closed surface.
In general, however, there may be no such path, as McDuff showed in~\cite{87}.

Take $M = S^2 \times S^2 \times T^2$ with the standard split symplectic form~$\omega$
giving all factors area one.
Define 
$$
\Psi \bigl( z, w, (s,t) \bigr)  \,=\, \bigl( z, R_{z,t}(w), s,t  \bigr)
$$
where $R_{z,t}$ is the rotation by angle~$2\pi t$ of the round sphere
about the axis through~$z,-z$.
{\it Then the symplectic forms $\omega_k = (\Psi^k)^* \omega$, $k \in \ZZ_{\geq 0}$,
are all cohomologous, and they can be joined by a path of symplectic forms, but by no such path in
the same cohomology class.} 
The last part of the statement was McDuff's first result 
obtained by using $J$-holomorphic curves.
She studied the space of $\omega_k$-tame $J$-holomorphic spheres 
in the class of the first $S^2$-factor,  
and showed that such spheres wrap $k$-times around the second factor.
Since a diffeomorphism isotopic to the identity acts trivially on homology, 
an isotopy deforming $\omega_k$ to $\omega_\ell$ can therefore only exist if $k = \ell$.

Also in \cite{87} McDuff improved this construction to obtain cohomologous symplectic forms
on an 8-dimensional symplectic manifold that are not even diffeomorphic.

\subsection*{Disconnected contact boundaries}

The natural boundaries of symplectic manifolds are those of contact type, 
meaning that there exists a vector field~$X$ defined near~$\pp M$ that is transverse to~$\pp M$, 
pointing outwards, and is conformally symplectic: $\cl_X \omega = d \iota_X \omega = \omega$.
Symplectic manifolds with boundary of contact type are analogous in many ways to
complex manifolds with pseudo-convex boundaries, as was shown by Eliashberg and Gromov~\cite{El89, EG89, El90}.
In the latter situation, the boundary is always connected. 
However, McDuff~\cite{91b} explicitly constructed a compact symplectic 4-manifold 
whose boundary is of contact type and disconnected.

She starts with the cotangent bundle $T^*\Sigma$ over a closed orientable surface of genus~$\geq 2$,
endowed with its canonical symplectic form~$d \lambda$, where $\lambda = \sum_i p_i \1 dq_i$.
Also take a Riemannian metric of constant curvature on~$\Sigma$. 
It induced the connection $1$-form~$\beta$ of the Levi-Civita connection on~$T^*\Sigma$
and the radial coordinate~$r$ on the fibers. 
Then one finds smooth functions $f,g$ on $[0,\infty)$ such that 
on the annulus bundle $\{ x \in T^*\Sigma \mid \frac 12 \leq r(x) \leq 1\}$ the 2-form 
$$
\omega \,=\, d \bigl( f(r) \2 \beta + g(r) \1 \lambda \bigr)
$$
is symplectic and makes the boundary of contact type.

\subsection{The structure of rational and ruled symplectic 4-manifolds}
\label{ss:structure}
The classification of compact complex surfaces is an old and beautiful topic
in complex and K\"ahler geometry. 
Before 1990, the only result on {\it symplectic} 4-manifolds in this direction
was a theorem of Gromov for the complex projective plane.
In \cite{90}, McDuff proved the following generalization:
{\it If a closed symplectic 4-manifold contains a symplectically embedded sphere
with non-negative self-intersection number, 
then it is symplectomorphic to either $\CP^2$ with its standard K\"ahler structure, 
to a ruled symplectic manifold, 
or to a symplectic blow-up of one of these manifolds.}
Here, a ruled symplectic 4-manifold is the total space of an $S^2$-fibration 
over a closed oriented surface with 
a symplectic structure that is non-degenerate on the fibers.

In later work with Lalonde~\cite{LaMc96a}, McDuff classified ruled symplectic surfaces:
{\it If $(M, \omega)$ is a ruled symplectic 4-manifold, then $\omega$ is determined up to symplectomorphism 
by its cohomology class and is isotopic to a standard K\"ahler form on~$M$.}
Li--Liu~\cite{LiLiu95} complemented this result by showing that if $(M, \omega)$ is the total space of an 
$S^2$-fibration over a closed surface, then there is a ruling of~$M$ by symplectic spheres, i.e.\ 
$(M, \omega)$ is a ruled symplectic manifold.
See~\cite{LaMc96b} for a survey on this classification.

An elementary but important point in the proofs is that two cohomologous symplectic forms
$\omega_0$ and~$\omega_1$ are diffeomorphic  
if they tame the same almost complex structure~$J$.
Indeed, all the cohomologous forms $\omega_t = (1-t) \omega_0 + t \omega_1$, $t \in [0,1]$, 
then tame~$J$ and hence are symplectic, 
and therefore $\omega_0$ and~$\omega_1$ are diffeomorphic by Moser's argument.

These works showed that symplectic geometry has something to say about 4-manifolds, 
and they helped establish symplectic geometry as one of the core geometries.

\m
McDuff has also done much interesting work on the topology of 
the group of symplectomorphisms for several symplectic manifolds. 
For $S^2$-fibrations over~$S^2$ with any symplectic form
she has in particular shown in joint work with Abreu~\cite{AbMc00} 
that two symplectomorphisms are isotopic through symplectomorphisms
whenever they are isotopic through diffeomorphisms.

\subsection{Hofer geometry} 
\label{ss:Hofer}

Recall that any symplectic manifold $(M,\omega)$ locally looks like 
the standard symplectic vector space of the same dimension.
Our dear geometric intuition from everyday life, 
so useful in Riemannian geometry to see distances and curvatures, 
is thus useless in symplectic geometry. 
However, on the automorphism group $\Hamc (M,\omega)$ of Hamiltonian diffeomorphisms
that are generated by compactly supported functions, 
there is a bi-invariant Finsler metric, that to some extent serves as a substitute
for the absence of local geometry in $(M,\omega)$.
Given a time-dependent Hamiltonian function $H \colon M \times [0,1] \to \RR$
with compact support,
take the integrated oscillation
$$
\| H \| \,:=\, \int_0^1 \left( \max_{x \in M} H(x,t) - \min_{x \in M} H(x,t) \right) dt .
$$
For $\gf \in  \Hamc (M,\omega)$ now define 
$$
d(\gf, \id) \,=\, \inf_H \|H\|
$$
where $H$ runs over all compactly supported Hamiltonian functions 
whose time-1 flow map is~$\gf$.
Then $d(\gf,\psi) = d(\gf \psi^{-1},\id)$ defines a bi-invariant metric on $\Hamc (M,\omega)$.
The only difficult point to check is non-degeneracy. 
This was done by Hofer~\cite{Ho90} for~$\RR^{2n}$ by variational methods, 
by Polterovich~\cite{Pol93} for tame rational symplectic manifolds 
by using a rigidity property of Lagrangian submanifolds,
and for all symplectic manifolds by Lalonde--McDuff~\cite{LaMc95} who used
their General Nonsqueezing Theorem discussed in~\S \ref{ss:nonsq} 
and the symplectic folding construction:

\begin{figure}[h] 
 \begin{center}
     \psfrag{a}{$a$}
     \psfrag{a2}{$\frac a2$}
     \psfrag{M}{$M$}
     \psfrag{R}{$\RR^2$}
     \psfrag{s}{$\sigma \times \id$}
     \psfrag{l}{$\lambda$} 
     \psfrag{t}{$\tau \times \id$} 
  \leavevmode\epsfbox{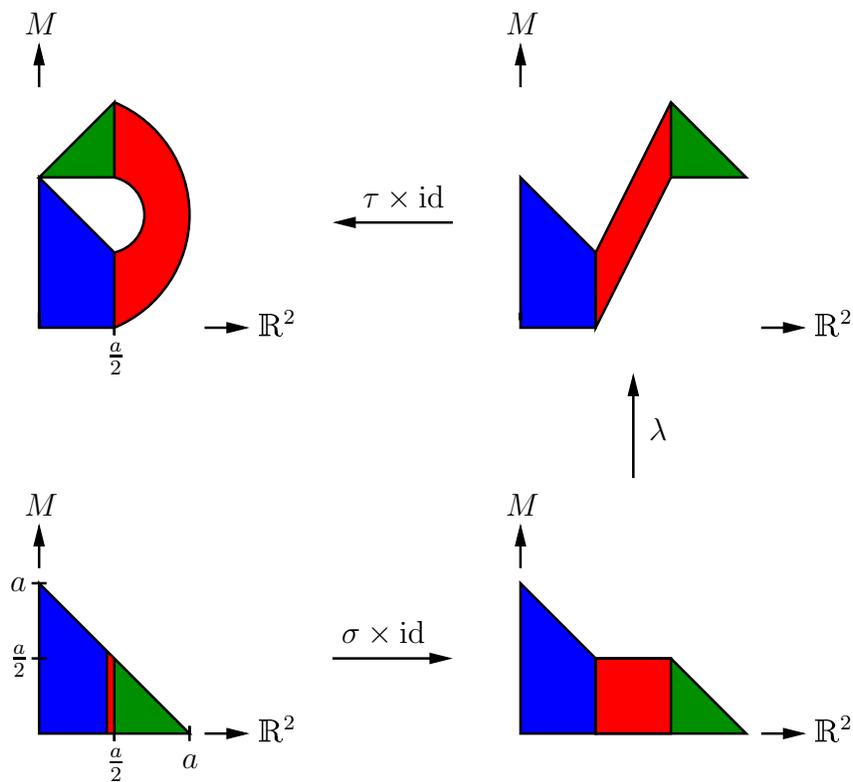}  
 \end{center}
 \caption{The symplectic folding construction, schematically}
 \label{fig.folding}
\end{figure}
%

\ni
Assume that $\gf \neq \id$. 
Then we find a symplectically embedded ball $\B^{2n}(a)$ in~$M$
such that $\gf (\B^{2n}(a)) \cap \B^{2n}(a) = \emptyset$.
Now consider the ball 
$$
\B^{2n+2}(a) \,\subset\, \B^2(a) \times \B^{2n}(a) \,\subset\, \RR^2 \times M ,
$$
and assume that $\gf$ is generated in time~$1$ by~$H$.
In three steps, the small green ball $\B^{2n+2}(\frac a2) \subset \B^{2n+2}(a)$ is folded on top of 
its complement,
as illustrated in Figure~\ref{fig.folding}: 
First, one views $\B^{2n+2}(a)$ as a $\B^{2n}$-fibration over the disc~$\B^2(a)$
and separates the small fibres from the large ones.
In the key step, one then uses the flow $\phi^t_H$, $t \in [0,1]$,
to lift the small green ball. 
The projection of the red band in the image of~$\lambda$ to~$\RR^2$
has area~$\|H\|$. 
Since  $\gf (\B^{2n}(\frac a2))$ is disjoint from $\B^{2n}(a)$,
one can thus turn the green ball over the blue part to obtain an embedding
$$
\B^{2n+2}(a) \,\se\, \B^2 \bigl( \tfrac a2+\|H\| \bigr) \times M .
$$
Hence $\|H\| \geq \frac a2$ by the General Nonsqueezing Theorem.
You can readily grasp the construction from the three figures on pages~473, 474, 475 of~\cite{McSal.book1}.
Symplectic folding found many other applications to symplectic geometry, 
see for instance~\cite{Hi15}.

In \cite{LaMc95b}, Lalonde and McDuff made a deep study of geodesics in this Finsler geometry 
on~$\Hamc (M,\omega)$,
that lead to completely new geometric intuitions on this group. 
It then became possible to think about Hamiltonian dynamics 
in geometric terms such as geodesic, conjugate point, cut-locus, etc. 
See Polterovich's book~\cite{Pol01} for more on this.

\subsection*{The books.}
While in her papers Dusa shows herself to be an impressive and creative problem solver,
the two books she wrote with Dietmar Salamon were (and remain) 
crucial for the foundation of symplectic geometry and its dissemination.
The ``Introduction to symplectic topology''~\cite{McSal.book1}
explains the methods and results of the field in clear and modern 
geometric language.
A key for the success of this book is that it discusses the main results in the most important 
and typical cases, without striving for generality.
I was very lucky that this book came out (in 1995) just when I wanted to learn the subject.
It is one of the main reasons for the transformation of the then 
small community of symplectic geometers into a large family.
In the third edition from~2017, the authors added a chapter with 54~open problems,
proving that symplectic geometry is not a closed chapter but
rather remains an exploding field.

The second book \cite{McSal.book2} provides a rigorous foundation of the theory of
$J$-holomorphic curves and explains their applications to symplectic topology.
The exposition is so precise and to the point that one can often just cite the
result one needs.
This book transformed the formerly somewhat romantic theory of $J$-holomorphic curves 
into a well established tool of enormous impact.

\b
In addition to her incredible mathematical legacy, 
Dusa's extraordinary generosity to young (and not so young) researchers, 
her enthusiasm and her (sometimes overwhelming) energy, 
and her heartfelt commitment in scientific, political, and social issues 
should be recognized.

\b \ni
{\bf Acknowledgments.}
I am grateful to the two referees for the improvements, 
and to Jesse Litman for her patient help with the English.

\end{document}